\documentclass[smallcondensed]{svjour3} 

\smartqed

\usepackage[bookmarks=true]{hyperref}  % use this for the draft

\usepackage{ifpdf}
\ifpdf
   \usepackage[pdftex]{graphicx}
   \usepackage{epstopdf}
\else
   \usepackage{graphicx}
\fi

\usepackage{tikz}

\usepackage{color}
\usepackage{amssymb,amsmath}

\begin{document}

\title{Richardson Extrapolation for Linearly Degenerate Discontinuities\thanks{This work was performed under the auspices of the
  U.S. Department of Energy by Lawrence Livermore National Laboratory
  under Contract DE-AC52-07NA27344 and was funded by 
  the Uncertainty Quantification Strategic Initiative Laboratory Directed 
  Research and Development Project at LLNL under
  project tracking code 10-SI-013, by 
  DOE contracts from the ASCR Applied Math Program, and by Los Alamos National Laboratory under Contract DE-AC52-06NA25396.}}

\author{J. W. Banks \and T. D. Aslam}

 \institute{J. W. Banks 
   \at Center for Applied Scientific Computing, L-422 \\
   Lawrence Livermore National Laboratory\\
   Livermore, CA 94551\\
   Tel.: 925-423-2697 \\
   \email{banks20@llnl.gov}
      \and T. D. Aslam
    \at Weapons Experiments, MS P952 \\
    Los Alamos National Laboratory\\
    Los Alamos, NM 87545\\
    Tel.: 505-667-1367\\
    \email{aslam@lanl.gov}}

\maketitle

\begin{abstract}
In this paper we investigate the use of Richardson extrapolation to estimate the convergence rates for numerical solutions to wave propagation problems involving discontinuities. For many cases, we find that the computed results do not agree with the a-priori estimate of the convergence rate. Furthermore, the estimated convergence rate is found to depend on the specific details of how Richardson extrapolation was applied, in particular the order of comparisons between three approximate solutions can have a significant impact. Modified equations are used to analyze the situation. We elucidated, for the first time, the cause of apparently unpredictable estimated convergence rates from Richardson extrapolation in the presence of discontinuities. Furthermore, we ascertain one correct structure of Richardson extrapolation that can be used to obtain predictable estimates. We demonstrate these results using a number of numerical examples.

% \PACS{PACS code1 \and PACS code2 \and more}
%\subclass{65M12 \and 65M08 \and 65M06}
\end{abstract}

%\begin{AMS}
%65M12, 65M08, 65M06
%\end{AMS}

\newcommand{\cfl}{\lambda}
\newcommand{\erf}{\operatorname{erf}}
\newcommand{\erfc}{\operatorname{erfc}}
\newcommand{\hu}{\hat{u}(x,t)}
\newcommand{\bu}{u(x,t)}
\newcommand{\bU}{U(x,t)}

\newcommand{\ha}{h_1}
\newcommand{\hb}{h_2}
\newcommand{\hc}{h_3}
\newcommand{\norm}[1]{{\left|\left| #1 \right|\right|}}
\newcommand{\RE}[3]{{{\mathcal{R}}\left(#1,#2,#3\right)}}

\newcommand{\junk}[1]{{}}

\newcommand\mstrut{\rule[-1.5ex]{0pt}{4ex}}
\newcommand\bstrut{\rule[-1.2ex]{0pt}{0pt}}

\newcommand{\utt}{\frac{\partial^2 u}{\partial t^2}}
\newcommand{\uxx}{\frac{\partial^2 u}{\partial x^2}}
\newcommand{\uyy}{\frac{\partial^2 u}{\partial y^2}}

\newcommand{\uy}{\frac{\partial u}{\partial y}}
\newcommand{\ut}{\frac{\partial u}{\partial t}}

\newcommand{\xv}{{\bf x}}

\newcommand{\Rich}{Richardson extrapolation }

%-------------- Introduction
\section{Introduction}
\label{sec:intro}
Estimating the error in numerical approximations to solutions of partial differential equations is important for many reasons. In order to be useful, numerical simulations should be accurate in some measurable norm. A particular application could require the absolute error to be less than a certain level. Other applications might use an estimate of the numerical error to help guide decisions about the most cost effective way to spend scarce resources, for example to choose between higher resolution simulations, to include more physical processes into the model, or to produce more samples for a statistical analysis. Still other applications may use estimates of the error directly for uncertainty quantification purposes. There are many approaches to error estimation found in the literature. Intrusive techniques such as adjoint error estimators~\cite{giles02,estep00}, error transport~\cite{hay06,banks12b}, or finite element residual and recovery methods~\cite{ainsworth97} are extremely powerful. However, because they are intrusive they require access and modification of the source code. This is often not possible for theoretical, practical, or sometimes legal reasons. 

Non-intrusive techniques are those that require only the ability to produce multiple simulation results, but do not require modifications of the source code. Error estimation through Richardson extrapolation is one commonly used non-intrusive technique and essentially relies on asymptotic properties of numerical approximation. Asymptotically correct in this context refers to the fact that expending a certain amount more computational effort yields a predictable increase in the accuracy of the result. Richardson extrapolation estimates can be based on varying the order of approximation, varying the resolution of the grid, or a combination of both. For many applications, \Rich has been shown to yield very good results. However, the behavior of \Rich error estimates for simulations of solutions with discontinuities is known to be problematic~\cite{roy10}. Discontinuous solutions are common in many physical systems, e.g. near material interfaces in electromagnetics or shock and contact waves for fluid mechanics, and so understanding the behavior or \Rich in theses cases can be important. There have been attempts to introduce richer ansatz when using \Rich in the presence of discontinuities, and these have yielded varying degrees of success. However, there has been very little progress on understanding the fundamental sensitivity of \Rich error estimates in the presence of jumps. The current work is a step toward filling that gap. 

In this paper, we investigate one particular realization of Richardson extrapolation error estimation for linear wave propagation with discontinuous solutions. Linear jumps are important in their own rite, for example wave propagation through solids~\cite{niethammer95,appelo12}, and for their significant role in nonlinear problems, for example contact and slip surfaces in the Euler equations~\cite{lax72,banks08b}. In this manuscript, we build on the previous work of~\cite{banks08a} which analyzed convergence rates for approximate solutions of linear advection where the exact solution contained jumps. That work used modified equations to argue that the expected rate of convergence for a nominally $p^{\hbox{\footnotesize{th}}}$ order method in the presence of a linear jump discontinuity is $p/(p+1)$. In the current work we use the structure of the modified equation solutions to discuss the expected behavior of \Rich error estimates. We show that under certain conditions one can expect to obtain the $p/(p+1)$ rate. In addition, we show why the method can fail to obtain the correct result if these conditions are not met.

The remainder of this paper is structured as follows. Section~\ref{sec:reminder} discusses some preliminaries, and provides a very brief overview of the Richardson extrapolation technique for error estimation. A motivating example problem demonstrating the difficulty associated with using \Rich error estimation for linear jumps is presented in section~\ref{sec:sosup}. Section~\ref{sec:model} gives an analysis of a simplified model problem of 1D linear advection. In Section~\ref{sec:firstOrder} we apply \Rich to approximations generated by a first-order upwind method, and discuss the results. The technique is found to be effective for this case, and an analysis explaining this surprising result is presented. That analysis is extended in Section~\ref{sec:HO} to discuss the case of high-order linear schemes. This analysis reveals that one particular instantiation of \Rich produces the expected convergence rate while others may not. Section~\ref{sec:examples} demonstrates the theory for upwind discretizations of order $2$, $4$, and $6$, as well as the case of a high-resolution nonlinear TVD discretization. Additional details of the inner workings are presented for the second-order case. In Section~\ref{sec:tophat}, the motivating example of Section~\ref{sec:sosup} is revisited, and the conclusions from the 1D analysis are shown to hold even for this more complex 2D case. Concluding remarks are presented in Section~\ref{sec:conclusions}.

%-------------- RE
\section{Preliminaries and Richardson extrapolation for smooth problems}
\label{sec:reminder}
Richardson extrapolation is a commonly used technique for error estimation, and many variations exist. For a good overview of the technique refer to~\cite{roy10}. Here we focus on one particular approach to \Rich which uses numerical approximations at three grid resolutions obtained using the same numerical technique.  Even within this particular flavor of the approach, there are essentially three possible realizations. In this section we review the approach and present the three choices.  In what follows, we consider numerical approximations to the solution of a partial differential equation (PDE). Boundary conditions are an important aspect of many numerical simulations, but are not critical to the present discussion. Thus consider the spatial domain $\xv\in\Omega$, and introduce a spatial discretization with uniform grid spacing $h$. 

Consider a set of numerical approximations given by $u_{h_M}(\xv,t)\approx u_e(\xv,t)$ where $u_e(\xv,t)$ is the exact solution, and $h_M$ indicates the size of the mesh. We consider performing an estimate at some time $t=t_f$, and whenever the time argument is not included, it is assumed to imply $t=t_f$ (i.e. $u(\xv) = u(\xv,t_f)$). Let an estimated convergence rate be denoted by $\RE{u_{h_1}}{u_{h_2}}{u_{h_3}}$ where the various $u_{h_M}$ are numerical approximations obtained using grid spacing $h_M$, and $\RE{u_{h_1}}{u_{h_2}}{u_{h_3}}=\sigma$ is the solution of the scalar equation $f(\sigma;u_{h_1},u_{h_2},u_{h_3})=0$, where 
\begin{equation}
  f(\sigma;u_{h_1},u_{h_2},u_{h_3}) =  \frac{\norm{u_{h_1}(\xv)-u_{h_2}(\xv)}}{\norm{u_{h_2}(\xv)-u_{h_3}(\xv)}}-\frac{\left| h_1^\sigma-h_2^\sigma \right|}{\left| h_2^\sigma-h_3^\sigma\right|}.
  \label{eq:RE_def}
\end{equation}
For the purposed of the remainder of this paper we will assume that $\norm{\cdot}$ indicates a discrete approximation to the $L_1$ norm. This is the norm which is most often considered when discussing hyperbolic equations with discontinuities. For a given set of three numerical approximations, there are essentially three distinct ways that the estimate can be computed $\RE{u_{h_1}}{u_{h_3}}{u_{h_2}}$, $\RE{u_{h_1}}{u_{h_2}}{u_{h_3}}$, and $\RE{u_{h_2}}{u_{h_1}}{u_{h_3}}$. As  shown below, this distinction is irrelevant for smooth problems. However, it will become important for solutions with discontinuities. For smooth problems, the assumption underlying basic Richardson extrapolation error estimation is that a given numerical approximation $u_{h_M}(\xv,t)$ is related to the exact solution $u_e(\xv,t)$ as
\[
  u_{h_M}(\xv,t) = u_e(\xv,t)+c(\xv,t)h_M^p +\ldots
\]
where $p$ is the formal order of accuracy of the approximation, and $c(\xv,t)$ is an order one function which is independent of the mesh parameters.  High-order terms are ignored, and the difference between approximations at two resolutions is
\begin{align*}
  u_{h_1}(\xv)-u_{h_2}(\xv) & = u_e(\xv)+c(\xv)h_1^p-u_e(\xv)-c(\xv)h_2^p \\
  & = c(\xv)\left(h_1^p-h_2^p\right).
\end{align*}
For any three resolutions then we find that 
\begin{align}
  \frac{\norm{u_{h_1}(\xv)-u_{h_2}(\xv)}}{\norm{u_{h_2}(\xv)-u_{h_3}(\xv)}} & = \frac{\norm{c(\xv) \left(h_1^{p}-h_2^{p}\right)}}{\norm{c(\xv) \left(h_2^{p}-h_3^{p}\right)}}  \nonumber \\
  & = \frac{\norm{c(\xv)}\, \lvert h_1^{p} -h_2^{p} \rvert}{\norm{c(\xv)}\, \lvert h_2^{p} -h_3^{p} \rvert} \nonumber \\
  & = \frac{\lvert h_1^{p} -h_2^{p}\lvert}{\lvert h_2^{p} -h_3^{p}\lvert}. \label{eq:smooth}
\end{align}
Assuming $h_1\ne h_2 \ne h_3$, Equation (\ref{eq:smooth}) and its counterparts can easily be used to show that $\RE{u_{h_1}}{u_{h_3}}{u_{h_2}}=\RE{u_{h_1}}{u_{h_2}}{u_{h_3}}=\RE{u_{h_2}}{u_{h_1}}{u_{h_3}}=p$. Note that in principle one can then use the computed convergence rate $\sigma$ in order to estimate the exact solution and obtain field estimates of the error. Such an approach is presented in detail in~\cite{henshaw06} and~\cite{banks09a}. Also note that there is the possibility multiple roots in (\ref{eq:RE_def}), but this situation is easily recognized in practice and so we do not discuss this further.

%-------------- sosup
\section{A motivating example}\label{sec:sosup}
In order to motivate the need for detailed understanding of Richardson extrapolation error estimators for linear wave propagation problems involving discontinuities, we revisit the recently published article~\cite{banks12c}. In that paper, the authors develop a new class of numerical methods for second-order wave equations. In two space dimensions, schemes of order 1, 2, 4, and a nominally $2^{\hbox{nd}}$ order nonlinearly limited scheme are presented. For the linear schemes the order of accuracy is proved. In order to demonstrate the properties of the methods, a number of test problems are presented. Where the exact solution is known and sufficiently smooth, the theoretical accuracy is confirmed. This provides confidence that the computer code correctly implements the numerical methods. However, for arguably the most interesting test problem, the exact solution is not known and convergence was judged only visually. In fact, the reason Richardson extrapolation was not used to quantify convergence for that case was the inherent inconsistency produced by the method, and the subsequent difficulty in interpreting the results. That difficulty is reproduced here.

Consider the initial boundary value problem
\begin{align}
   &\utt  =  \frac{1}{4} \uxx + \uyy, \qquad (x,y) \in (-\pi,\pi) \times (-2\pi,0),  \label{eq:surfacePDEa} \\
   &u(x,y,0) = u_0(x,y), \quad \ut(x,y,0) = v_0(x,y), \\
   &\uy(x,0,t) = 0, \quad \uy(x,-2\pi,t)=0, \quad u(x-\pi,y,t) = u(x+\pi,y,t),                    \label{eq:surfacePDEc} 
\end{align}
where 
\[
  %\ut(x,y,0) = 0, \qquad 
  %u(x,y,0) = \left\{
  v_0(x,y) = 0, \qquad 
  u_0(x,y) = \left\{
  \begin{array}{lll}
    1 \quad & \text{ if $ x^2+(y+2)^2 < 1$},  \\
    0 \quad & \hbox{ otherwise. }
  \end{array}
  \right.
\]
\newcommand{\labelFont}{\small}
\begin{figure}[hbt]
\begin{center}
  \begin{tikzpicture}[scale=1]
  \useasboundingbox (0,5) rectangle (13,9.15);  % set the bounding box (so we have less surrounding white space)
    \draw(0,5) node[anchor=south west] {\includegraphics[height=3.9cm]{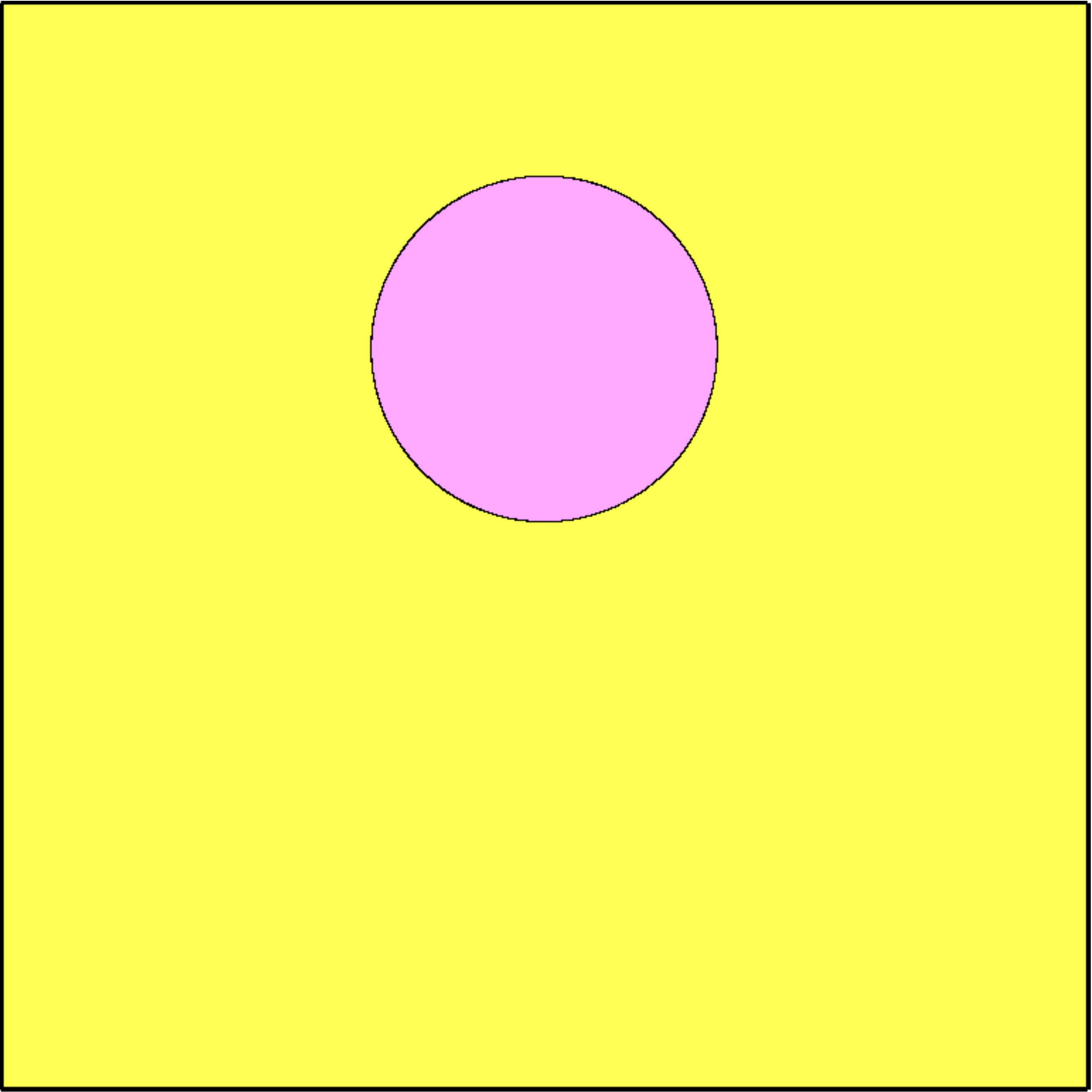}};
    \draw(4,5) node[anchor=south west] {\includegraphics[height=3.9cm]{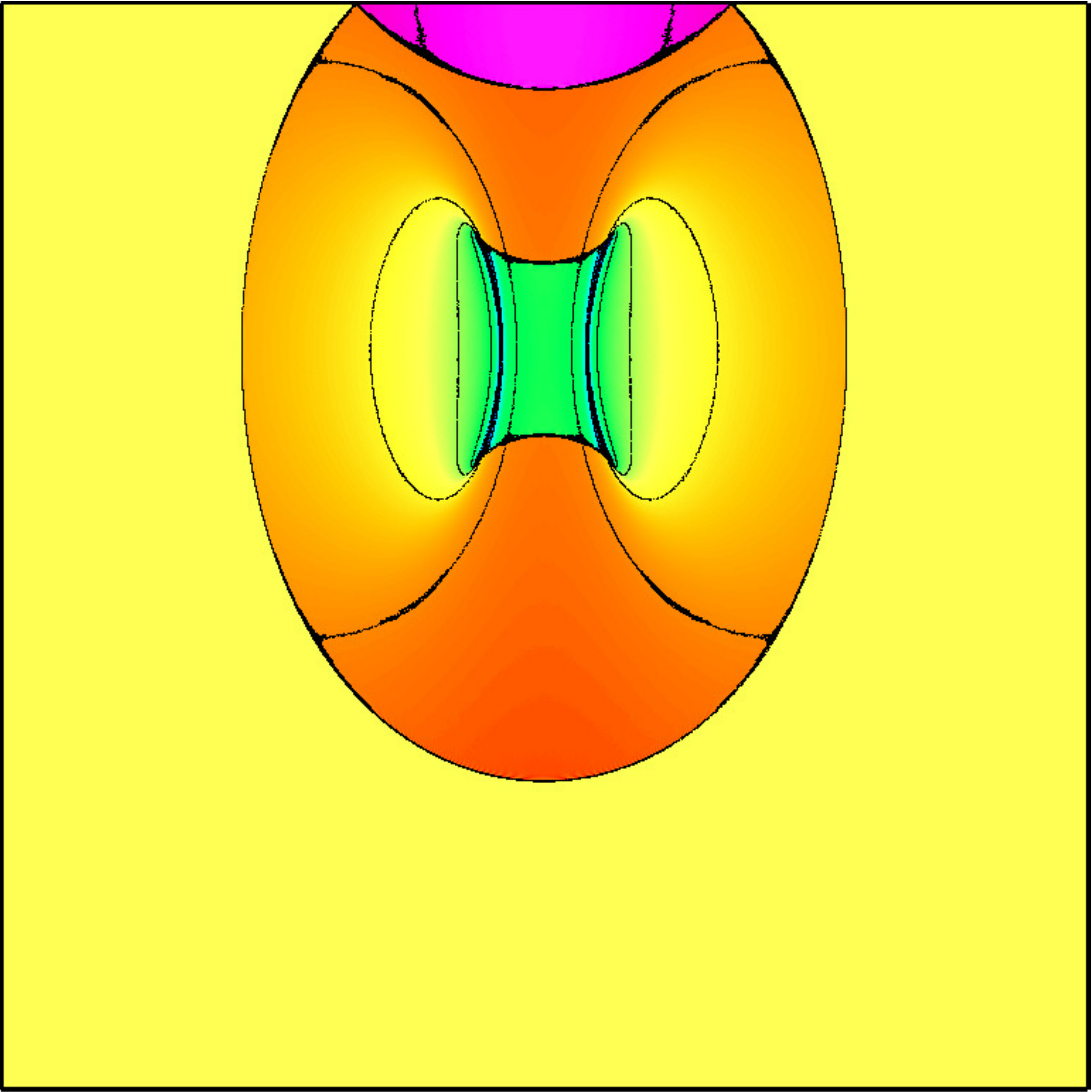}};
    \draw(8,5) node[anchor=south west] {\includegraphics[height=3.9cm]{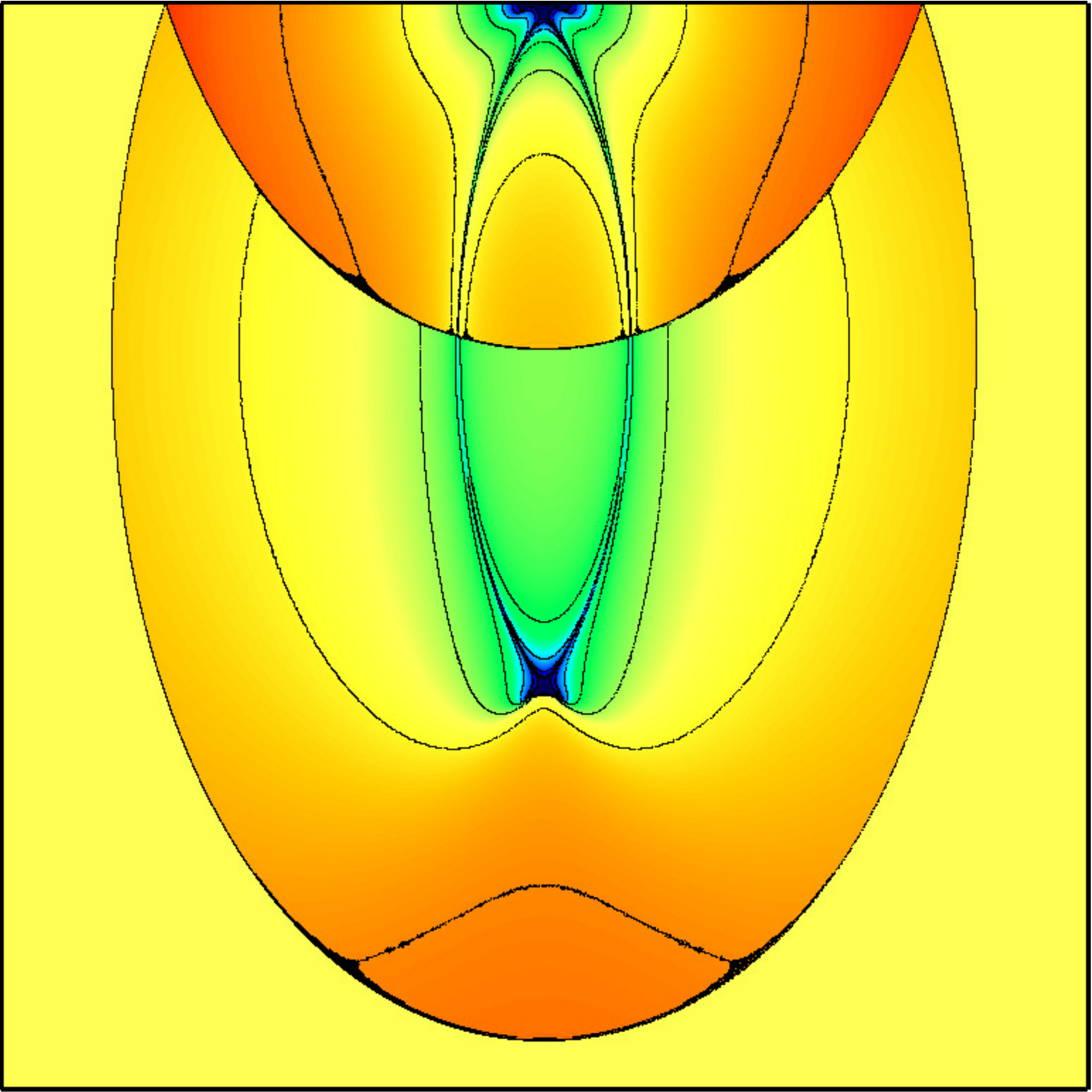}};
    \draw (12,5.05) node[anchor=south west,]  {\includegraphics[height=3.8cm]{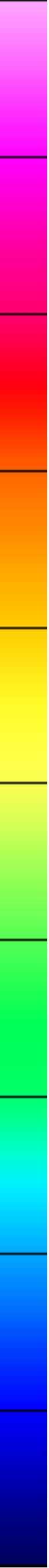}};
    \draw (12.5,5.25) node{\tiny$-1$};
    \draw (12.5,7.05) node{\labelFont$u$};
    \draw (12.5,8.85) node{\tiny $1$};
    \draw (.2,5.5) node[anchor=west, draw,fill=white, rounded corners,thick] {\scriptsize $t=0$};
    \draw (4.2,5.5) node[anchor=west, draw,fill=white, rounded corners,thick] {\scriptsize $t=1.5$};
    \draw (8.2,5.5) node[anchor=west, draw,fill=white, rounded corners,thick] {\scriptsize $t=3$};
    %
    % grid:
    %\draw[step=1cm,gray] (0,5) grid (13,9.25);
    %\draw[step=1cm,gray] (0,5) rectangle (13,9.15);
  \end{tikzpicture}
  \caption{Time evolution of the two-dimensional top-hat using a fourth-order accurate upwind method and $6400$ points in each coordinate direction. See~\cite{banks12c} for details concerning the method.} 
  \label{fig:topHat2D_evolution}
\end{center}
\end{figure}
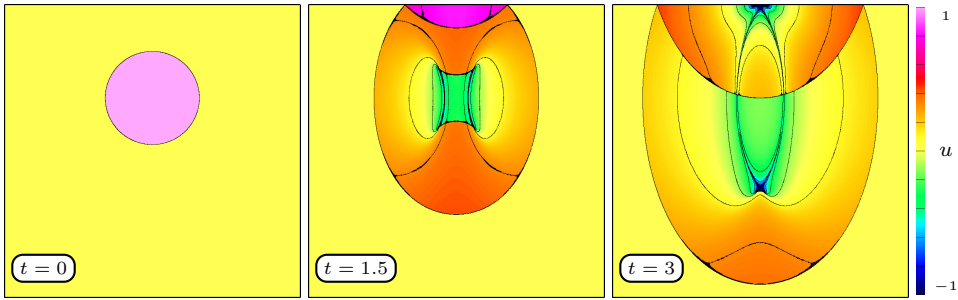
% ---------------------------------------------------------------------------------------------------------------

Figure~\ref{fig:topHat2D_evolution} shows the evolution of the numerical solution computed using the fourth-order accurate scheme from~\cite{banks12c} at times $t=0$, $t=1.5$, and $t=3.0$. The solution was computed using $6400$ points in the $x$- and $y$-directions. In order to judge convergence one could apply  Richardson extrapolation as discussed in Section~\ref{sec:reminder}.
\begin{table}[hbt]
\begin{center}
\begin{tabular}{|c||c|c|c|} \hline\hline
  \mstrut scheme & \footnotesize{$\RE{u_{h_1}}{u_{{h_1/2}}}{u_{h_1/4}}$} &\footnotesize{$\RE{u_{h_1}}{u_{h_1/4}}{u_{h_1/2}}$} & \footnotesize{$\RE{u_{h_1/2}}{u_{h_1}}{u_{h_1/4}}$} \\ \hline\hline
  \mstrut first order           & $0.44$ & $0.40$ & $0.50$ \\ \hline
  \mstrut second order    & $0.53$ & $-0.27$ & $1.92$ \\ \hline
  \mstrut high resolution & $0.62$ & $0.52$ & $0.79$ \\ \hline  
  \mstrut fourth order       & $0.64$ & $-0.50$ & $3.41$ \\ \hline\hline
\end{tabular}
\caption{Estimated convergence rates for the top-hat problem using upwind schemes and the three variants of Richardson extrapolation. The base resolution uses $400$ points in each coordinate direction, and uniform refinement is carried out using a ratio of $2$. Results are presented for numerical methods of order 1, 2, 4, and a high-resolution nonlinear variant of nominally second-order.}
\label{table:TH_naieve}
\end{center}
\end{table}
Table~\ref{table:TH_naieve} presents such a study using uniform refinement and shows results for the three variants of Richardson extrapolation as applied to approximations from numerical methods of order 1, 2, 4, and a high-resolution nonlinear variant of nominally second-order. As discussed in~\cite{banks08a,banks12c}, the expected convergence rate for this case is $p/(p+1)$, where $p$ is the nominal order of the method. The challenge in interpreting the results in Table~\ref{table:TH_naieve} is apparent as there is wide variations in the approximated convergence rate depending on which ordering within the Richardson extrapolation scheme is used. Indeed the fourth-order scheme yields results between $-.5$ and $3.41$, despite the fact that reasonably highly resolved simulations are being considered. Furthermore, we note that the results presented in the first column are in reasonable agreement with the expectation.

%-------------- model
\section{A model problem with discontinuity} 
\label{sec:model}
We now seek to understand the nature of Richardson extrapolation error estimation for problems with discontinuities, or other self similar behavior, using a simple model problem. Consider the one dimensional linear advection equation
\begin{equation}
 \frac{\partial}{\partial t} \bu+a \frac{\partial}{\partial x} \bu=0
 \label{eq:1Dadvection}
\end{equation}
with constant advection velocity $a>0$. A canonical model problem with discontinuity can be defined using the initial conditions 
\begin{equation}
  u(x,0)=\left\{
  \begin{array}{lcr}
    u_L \hspace{0.1in} & \hbox{for} & \hspace{0.1in} x < 0\\
    u_R \hspace{0.1in} & \hbox{for} & \hspace{0.1in} x \ge 0.\\
  \end{array}
  \right.
  \label{eq:jumpIC}
\end{equation}
The method of characteristics is used to define the exact solution for all $t>0$ as $u(x,t) = u(x-at,0)$, which applies also to discontinuous solution profiles using the notion of weak solutions~\cite{whitham74,lax72}. We note here that the expected order or accuracy for all numerical methods considered in this paper are confirmed using known smooth solutions.

%-------------- FOU
\subsection{First Order Upwind Discretization}
\label{sec:firstOrder}
Consider the first-order accurate explicit upwind scheme 
\begin{equation}
  u_i^{n+1}=u_i^n-\cfl \left[u_i^n-u_{i-1}^n\right]
  \label{eq:fou}
\end{equation}
where $u_i^n$ is a numerical approximation to $u(x_i,t^n)$ and the so called CFL number is $\cfl=\frac{a\Delta t}{h}$. The computational domain $[x_L,x_R]$ is a truncation of the infinite domain, and has been discretized with $x_i=x_L+i h$ where $h=(x_R-x_L)/(N-1)$ and $N$ an integer. Zero gradient conditions are applied at domain boundaries, but the nature of the specific problem being studied makes the details of these artificial boundary conditions unimportant. Time has been discretized as $t_n=n \Delta t$ with initial conditions $u(x,0)$ being given at $t=0$ with $u_L=-1$ and $u_R=1$.  Numerical stability is obtained for $\cfl\le1$. 

\subsubsection{Richardson extrapolation error estimation}
One can perform an estimate of the convergence rate using Richardson extrapolation and simply ignore the fact that the assumptions underlying the approach are not strictly valid for cases with discontinuities. We set $a=1$, choose a computational domain with $[x_L,x_R]=[-\pi,\pi]$, integrate to a time $t_f=2$, and use $\cfl=0.6$. Note that in order to ensure that the initial condition is applied consistently for all cases, the grid point located at $x=0$ is explicitly taken to have the initial value $u_L$.
%These choices are made to ensure that the discontinuity never lies on a cell boundary which can be problematic if finite-precision arithmetic leads to a jump in the initial discontinuity location during the refinement process. 
A series of approximations is generated using a uniform refinement process with a ratio $0\le r<1$ (i.e. $h_2=rh_1$ and $h_3=r^2h_1$), and starts with $51201$ points in the domain (i.e. $h_1=\frac{2\pi}{51200}$).
\begin{table}[hbt]
\begin{center}
\begin{tabular}{|c||c|c|c|} \hline\hline
  \mstrut$r$ & \footnotesize{$\RE{u_{h_1}}{u_{rh_1}}{u_{r^2h_1}}$} &\footnotesize{$\RE{u_{h_1}}{u_{r^2h_1}}{u_{rh_1}}$} & \footnotesize{$\RE{u_{rh_1}}{u_{h_1}}{u_{r^2h_1}}$} \\ \hline\hline
  \mstrut$\frac{1}{2}$   & $0.50$ & $0.50$ & $0.50$ \\ \hline
  \mstrut$\frac{2}{5}$   & $0.50$ & $0.50$ & $0.50$ \\ \hline
  \mstrut$\frac{1}{3}$   & $0.50$ & $0.50$ & $0.50$ \\ \hline
  \mstrut$\frac{2}{7}$   & $0.50$ & $0.50$ & $0.50$ \\ \hline
  \mstrut$\frac{1}{4}$   & $0.50$ & $0.50$ & $0.50$ \\ \hline\hline
\end{tabular}
\caption{Estimated convergence rates for the first-order upwind scheme for a solution with a discontinuity. The base resolution uses $51201$ points, and uniform refinement is carried out using a ratio of $r$. Results using the three independent variants of Richardson extrapolation are presented in the various columns.}
\label{table:1stOrderEstimates}
\end{center}
\end{table}
Table~\ref{table:1stOrderEstimates} shows the results using the three basic variants of Richardson extrapolation, and various choices of the uniform refinement ratio $r$. The table makes clear that the estimated convergence rate is $0.5$ for any choice, which agrees exactly with the expected convergence rate for a first-order scheme with a linear jump~\cite{banks08a}.

\subsubsection{Explanation of the result}
The fact that Richardson extrapolation seems to work, in terms of convergence rate estimates, for the first order upwind scheme even with a discontinuous solution is surprising. In order to understand this result we extend the analysis in~\cite{banks08a}. The approach makes use of modified equation for a more complete understanding of the behavior. The modified equation is a continuous PDE whose solution describes the approximate behavior of the well resolved components of the discrete solutions, and is derived by substituting continuous functions $U(x,t)$ into the discrete equation (\ref{eq:fou}) by setting $u_i^n=U(x_i,t^n)$, and expanding all terms in Taylor series about the point $(x,t)=(x_i,t^n)$. For the first-order upwind scheme the result is
\begin{equation}
  \frac{\partial}{\partial t}\bU+a \frac{\partial}{\partial x}\bU
    -\frac{a h}{2}\left(1-\cfl \right)\frac{\partial^2}{\partial x^2}\bU
    +\cdots=0.
  \label{eq:me_fou}
\end{equation}
Truncating Equation (\ref{eq:me_fou}) yields the advection-diffusion equation
\begin{equation}
  \frac{\partial}{\partial t}\bU+a \frac{\partial}{\partial x}\bU
    -\nu\frac{\partial^2}{\partial x^2}\bU=0
  \label{eq:modApprox}
\end{equation}
where $\nu = \frac{ah(1-\lambda)}{2}$. For discontinuous initial data (\ref{eq:jumpIC}), $U(x,0)=u_L$ for $x<0$ and $U(x,0)=u_R$ for $x  \ge 0$. The analytic solution to (\ref{eq:modApprox}) for $t>0$ is then found to be
\begin{equation}
  \bU=\frac{u_L+u_R}{2}+\frac{u_R-u_L}{2}\erf\left(\frac{x-a t}{\sqrt{4 \nu t}}\right)
  \label{eq:modExact}
\end{equation}
where $\erf(\zeta)$ is the error function
\[
  \erf(\zeta)=\frac{2}{\pi}\int^{\zeta}_0 e^{-\chi^2}d\chi.
\]
For additional details on this derivation refer to~\cite{banks08a}.

The analysis to follow assumes the use of the $L_1$ norm and sets $z=x-at_f$, $\delta_1=\sqrt{4\nu_1t_f}$ and $\delta_2=\sqrt{4\nu_2t_f}$. Furthermore, assume $h_1 > h_2$. As in~\cite{banks08a}, assume that the solution to the modified equation is an accurate approximation to the numerical solution so $u_h(x,t)=U(x,t)$. Following a similar line of reasoning as in Section~\ref{sec:reminder} gives
\begin{align*}
  \norm{u_{h_1}(x)-u_{h_2}(x)} & = \norm{\frac{u_L+u_R}{2}+\frac{u_R-u_L}{2}\erf\left(\frac{z}{\delta_1}\right)-\frac{u_L+u_R}{2}-\frac{u_R-u_L}{2}\erf\left(\frac{z}{\delta_2}\right)}\\
  & = \int_{-\infty}^{\infty} \left| \frac{u_R-u_L}{2}\left(\erf\left(\frac{z}{\delta_1}\right)-\erf\left(\frac{z}{\delta_2}\right) \right)\right|\, dz \\
  & = \left|u_R-u_L\right| \left( \int_0^{\infty} -\erf\left(\frac{z}{\delta_1}\right)+\erf\left(\frac{z}{\delta_2}\right)\, dz \right) \\
  & = \frac{2\sqrt{at_f}\left|u_R-u_L\right|}{\sqrt{\pi}}\left( \sqrt{\nu_1}-\sqrt{\nu_2} \right) \\
  & = \sqrt{\frac{at_f(1-\lambda)}{2\pi}} \left|u_R-u_L\right| \left(\sqrt{h_1}-\sqrt{h_2}\right).
\end{align*}
Therefore, under the assumption that the three numerical approximations have been obtained using the same CFL $\lambda$, the factor of $\sqrt{\frac{at_f(1-\lambda)}{2\pi}} \left|u_R-u_L\right|$ will appear in both the numerator and denominator when the ratio of the norms of differences is taken. As a result
\begin{align*}
   \frac{\norm{u_{h_1}(x)-u_{h_3}(x)}}{\norm{u_{h_2}(x)-u_{h_3}(x)}} 
  & = \frac{\left|\sqrt{h_1}-\sqrt{h_3}\right|}{\left|\sqrt{h_2}-\sqrt{h_3}\right|},
\end{align*}
and it is easy to verify that all three approaches to Richardson extrapolation will yield convergence at the expected rate of $0.5$.

\noindent{\em Remark:} Although the results presented in Table~\ref{table:1stOrderEstimates} use simulations with uniform refinement, this is not critical in the analysis for this case. In fact, it is primarily the monotone nature of the similarity solution which is responsible for the robust nature of the estimates. Uniform refinement was used in order to match the analysis for high-order schemes below, where uniform refinement is important.

%-------------- higher-order
\subsection{Arbitrary Order Linear Scheme}
\label{sec:HO}
In general, robust results like those represented in Table~\ref{table:1stOrderEstimates} are not expected. In this section we analyze why this is the case and determine a particular strategy which yields accurate estimates even in the presence of discontinuities (or other self-similar features). We restrict our attention to noncompressive stable $p^{\hbox{th}}$ order schemes for advection with modified equations of the form
\begin{equation}
  \frac{\partial}{\partial t}\bU+a \frac{\partial}{\partial x}\bU
    -\eta_h\frac{\partial^{p+1}}{\partial x^{p+1}}\bU
    +\cdots=0
    \label{eq:modHigher}
\end{equation}
where $\eta_h=\tilde{\eta}h^p$ and $\tilde{\eta}$ is a constant depending on the CFL $\cfl$. Following the analysis in~\cite{banks08a}, a simple change of variables is performed to translate into a frame of reference traveling with the wave 
\[
  \begin{array}{l}
    z=x-a t \smallskip \\
    \tau=t.
  \end{array}
\]
After dropping the higher order terms, (\ref{eq:modHigher}) becomes
\begin{equation}
  \frac{\partial}{\partial \tau}{U}(z,\tau)
    -\kappa_h\frac{\partial^{p+1}}{\partial z^{p+1}}{U}(z,\tau)=0
    \label{eq:modApproxHigher}
\end{equation}
where $\kappa_h$ is either plus or minus $\eta_h$ depending on the value of $p$. Similarity solutions can be sought with similarity variable
\begin{equation}
  \xi_h=\frac{z}{\sqrt[p+1]{\kappa_h t_f}}.
\end{equation}
We again assume that the solution of the modified equation is an accurate representation of the numerical approximation and set $u_h=U$. For jump initial condition (\ref{eq:jumpIC}) the solution can then be written in the general form
\begin{equation}
  u_h(\xi_h) = \frac{u_L+u_R}{2}+\frac{u_R-u_L}{2}S(\xi_h)
  \label{eq:genSimilarityForm}
\end{equation}
where $S$ is an approximation to the jump from $-1$ to $1$ (similar to an error function but perhaps with more complex  behavior). In general, $S$ will take the form of generalized hypergeometric functions which oscillate on one or both sides of the discontinuity. In a Richardson style error estimate, norms of the difference between two solutions will be used, and so we write
\begin{align*}
  \norm{u_{h_1}(x)-u_{h_2}(x)}
    & = \norm{\frac{u_L+u_R}{2}+\frac{u_R-u_L}{2}S(\xi_{h_1})-\frac{u_L+u_R}{2}-\frac{u_R-u_L}{2}S(\xi_{h_2})} \\
    & = \frac{\left|u_R-u_L\right|}{2}\int_{-\infty}^{\infty} \left| S\left(\frac{z}{\sqrt[p+1]{\kappa_{h_1} t_f}}\right)-S\left(\frac{z}{\sqrt[p+1]{\kappa_{h_2} t_f}}\right) \right| \, dz.
\end{align*}
Making the change of variables to
\begin{equation}
  \chi = \frac{z}{\sqrt[p+1]{\kappa_{h_1} t_f}}
  \label{eq:commonChange}
\end{equation}
gives
\begin{align}
  \norm{u_{h_1}(x)-u_{h_2}(x)}
    & = \left|u_R-u_L\right|\sqrt[p+1]{\kappa_{h_1}t_f}\int_{-\infty}^{\infty} \left| S\left(\chi\right)-S\left(\chi\sqrt[p+1]{\frac{h_1}{h_2}}\right) \right| \, d\chi.
    \label{eq:normDiff}
\end{align}
Now take the ratio of such norms to arrive at
\begin{align}
  \frac{\norm{u_{h_1}(x)-u_{h_2}(x)}}{\norm{u_{h_2}(x)-u_{h_3}(x)}} 
    & = \left(\frac{h_{1}}{h_{2}}\right)^{\frac{p}{p+1}}\frac{\int_{-\infty}^{\infty} \left| S\left(\chi\right)-S\left(\chi\sqrt[p+1]{\frac{h_1}{h_2}}\right) \right| \, d\chi}{\int_{-\infty}^{\infty} \left| S\left(\chi\right)-S\left(\chi\sqrt[p+1]{\frac{h_2}{h_3}}\right) \right| \, d\chi}.
  \label{eq:RE_new}
\end{align}
In general therefore, the Richardson estimate will depend on the ratio of integrals of scaled similarity functions
\begin{equation}
  \frac{\int_{-\infty}^{\infty} \left| S\left(\chi\right)-S\left(\chi\sqrt[p+1]{\frac{h_1}{h_2}}\right) \right| \, d\chi}{\int_{-\infty}^{\infty} \left| S\left(\chi\right)-S\left(\chi\sqrt[p+1]{\frac{h_2}{h_3}}\right) \right| \, d\chi}.
  \label{eq:primaryRatio}
\end{equation}
The function $S$ can be an extremely complex object, and so computing this ratio in closed form is in general impractical. In fact, the definition of $S$ is often found to make this ratio difficult to even estimate numerically due to ill-conditioning and finite precision arithmetic. However, for the case of uniform refinement when $h_3=rh_2=r^2h_1$, this ratio is simply unity. Therefore, for this special case, the estimated convergence rate will be given by $\RE{u_{h_1}}{u_{rh_1}}{u_{r^2h_1}}=p/(p+1)$. This is the expected convergence rate as discussed in~\cite{banks08a}. Notice that such cancelation requires both uniform refinement, and that the estimate be performed using differences of successive refinement. Other choices will in general not yield the rate $p/(p+1)$. 

To be clear, our analysis shows that for Richardson extrapolation of the form $\RE{u_{h_1}}{u_{rh_1}}{u_{r^2h_1}}$ where $r$ is a uniform refinement rate, the a-priori convergence rate of $p/(p+1)$ will be obtained even in the presence of discontinuities. Note here that $r$ is simply a uniform refinement, but its exact value is not important in the analysis. Other realizations of Richardson extrapolation may not yield this result, and the exact value of the computed estimate can vary widely from the expected rate.

%-------------- SOU
\subsection{Second Order Linear Scheme}
\label{subsec:Unlimited}
In order to demonstrate the implications of the analysis in Section~\ref{sec:HO}, consider the linear second-order upwind method
\begin{equation}
  u_i^{n+1}=u_i^n-\cfl \left[
    \left(u_i^n+\frac{1}{4}(1-\cfl)(u_{i+1}^n-u_{i-1}^n)\right)-
    \left(u_{i-1}^n+\frac{1}{4}(1-\cfl)(u_{i}^n-u_{i-2}^n)\right)
  \right].
  \label{eq:2ndUpwind}
\end{equation}
This is simply a second order unlimited Godunov method. As before, one can perform an estimate of the convergence rate using Richardson extrapolation. Again we set $a=1$, choose a computational domain with $[x_L,x_R]=[-\pi,\pi]$, integrate to a time $t_f=2$, and use $\cfl=0.6$. The series of approximations is generated using a uniform refinement process with a ratio $0\le r<1$ (i.e. $h_2=rh_1$ and $h_3=r^2h_1$), and starts with $51201$ points in the domain (i.e. $h_1=\frac{2\pi}{51200}$).
\begin{table}[hbt]
\begin{center}
\begin{tabular}{|c||c|c|c|} \hline\hline
  \mstrut$r$ & \footnotesize{$\RE{u_{h_1}}{u_{rh_1}}{u_{r^2h_1}}$} &\footnotesize{$\RE{u_{h_1}}{u_{r^2h_1}}{u_{rh_1}}$} & \footnotesize{$\RE{u_{rh_1}}{u_{h_1}}{u_{r^2h_1}}$} \\ \hline\hline
  \mstrut$\frac{1}{2}$   & $0.67$ & $0.14$ & $1.63$ \\ \hline
  \mstrut$\frac{2}{5}$   & $0.67$ & $0.22$ & $1.73$ \\ \hline
  \mstrut$\frac{1}{3}$   & $0.67$ & $0.34$ & $1.60$ \\ \hline
  \mstrut$\frac{2}{7}$   & $0.67$ & $0.41$ & $1.49$ \\ \hline
  \mstrut$\frac{1}{4}$   & $0.67$ & $0.47$ & $1.35$ \\ \hline\hline
\end{tabular}
\caption{Estimated convergence rates for the second-order upwind scheme for a solution with a discontinuity. The base resolution uses $51201$ points, and uniform refinement is carried out using a ratio of $r$. Results using the three independent variants of Richardson extrapolation are presented in the various columns.}
\label{table:2ndOrderEstimates}
\end{center}
\end{table}
Table~\ref{table:2ndOrderEstimates} shows the results using the three basic variants of Richardson extrapolation, and various choices of the uniform refinement ratio $r$. The table shows that the expected rate of $2/3$ is obtained as advertised when using the appropriate approach. Also shown in the table are the type of results that can be experienced if other approaches are used.

\subsubsection{A closer look at this case}
\begin{figure}
\begin{center}
  \includegraphics[width=3.0in]{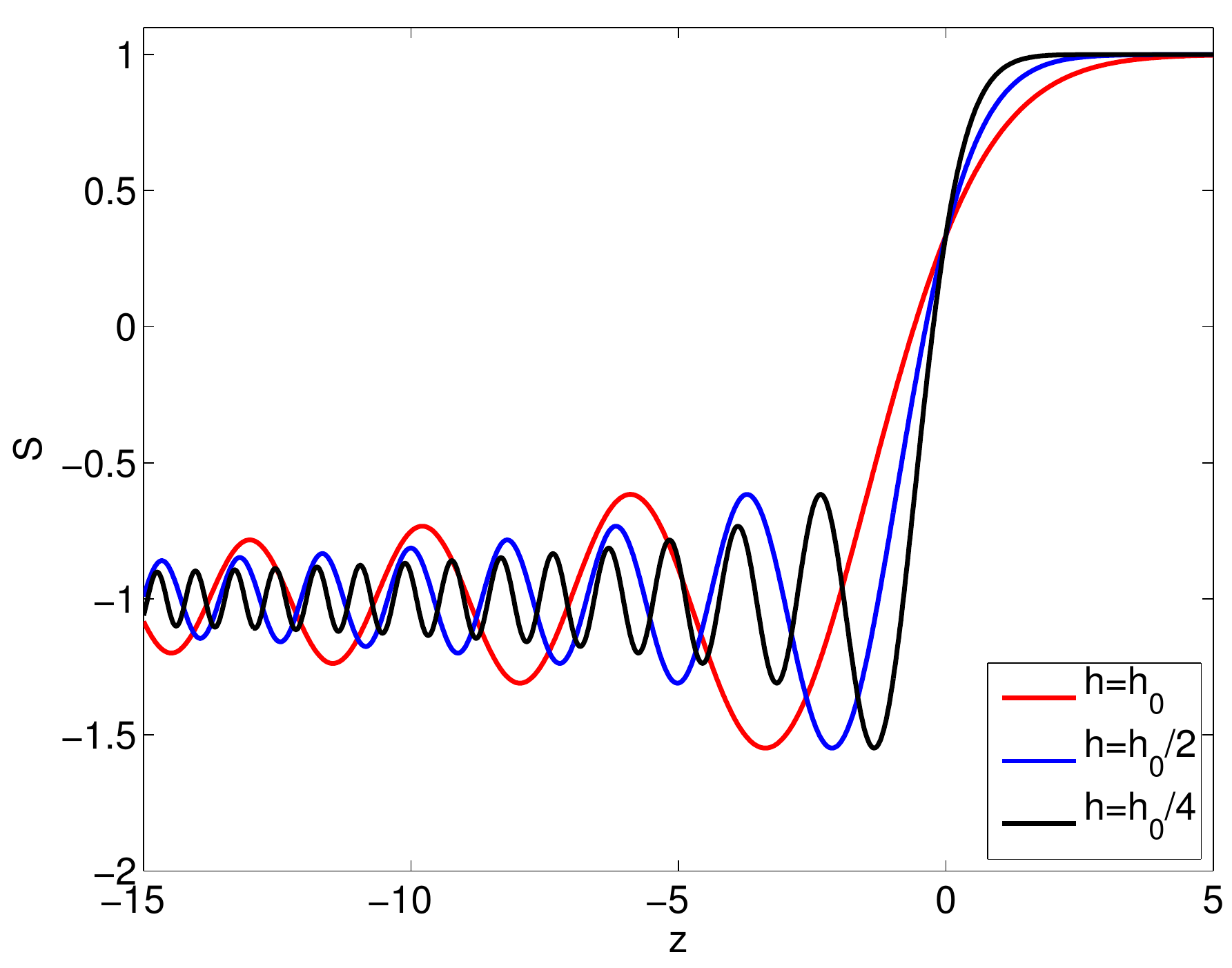}
  \caption{Similarity jumps for the unlimited second order scheme in the frame of reference  moving with the discontinuity.}
  \label{fig:souJump}
\end{center}
\end{figure}

\begin{figure}
\begin{center}
  \includegraphics[width=2.0in]{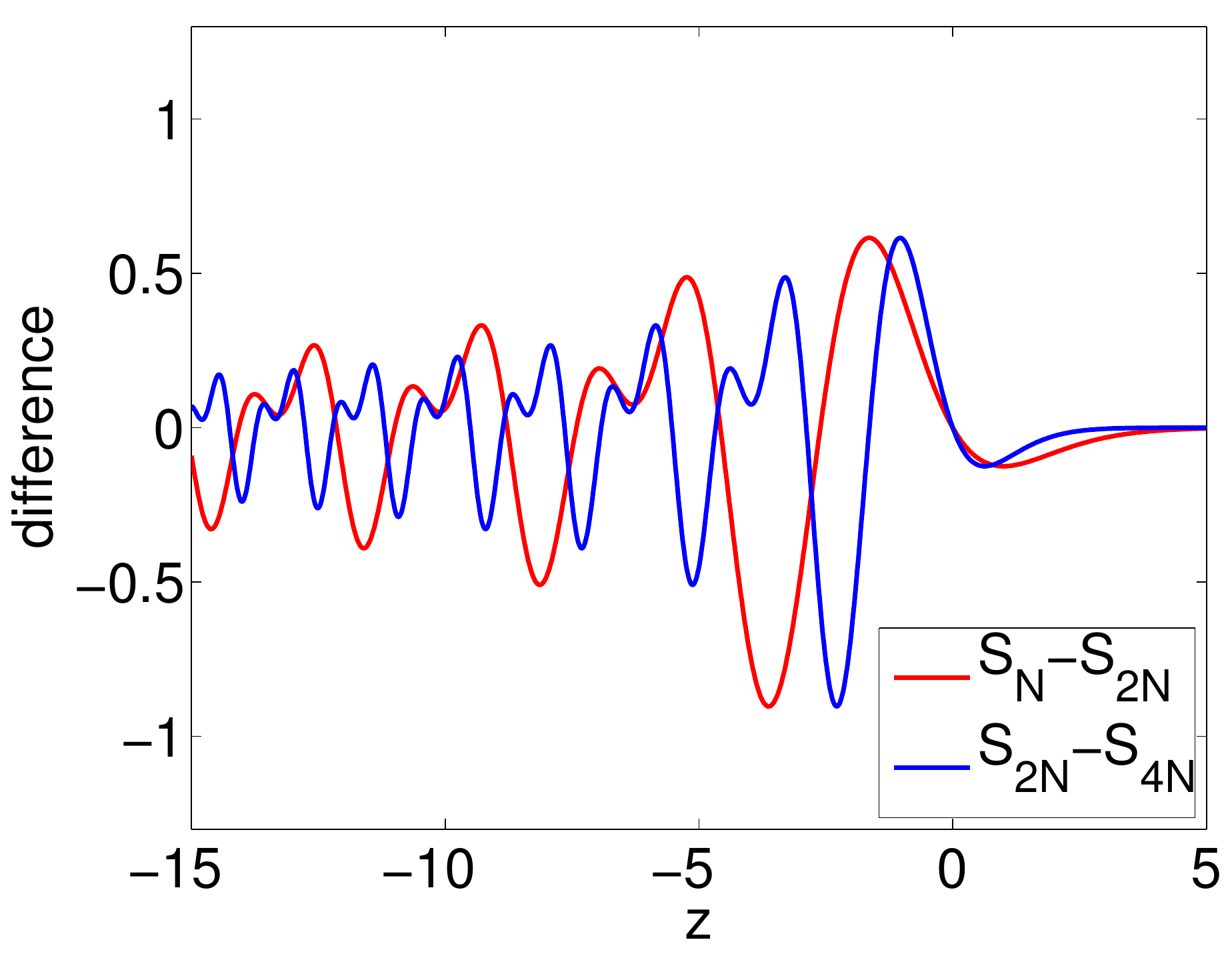}
  \includegraphics[width=2.0in]{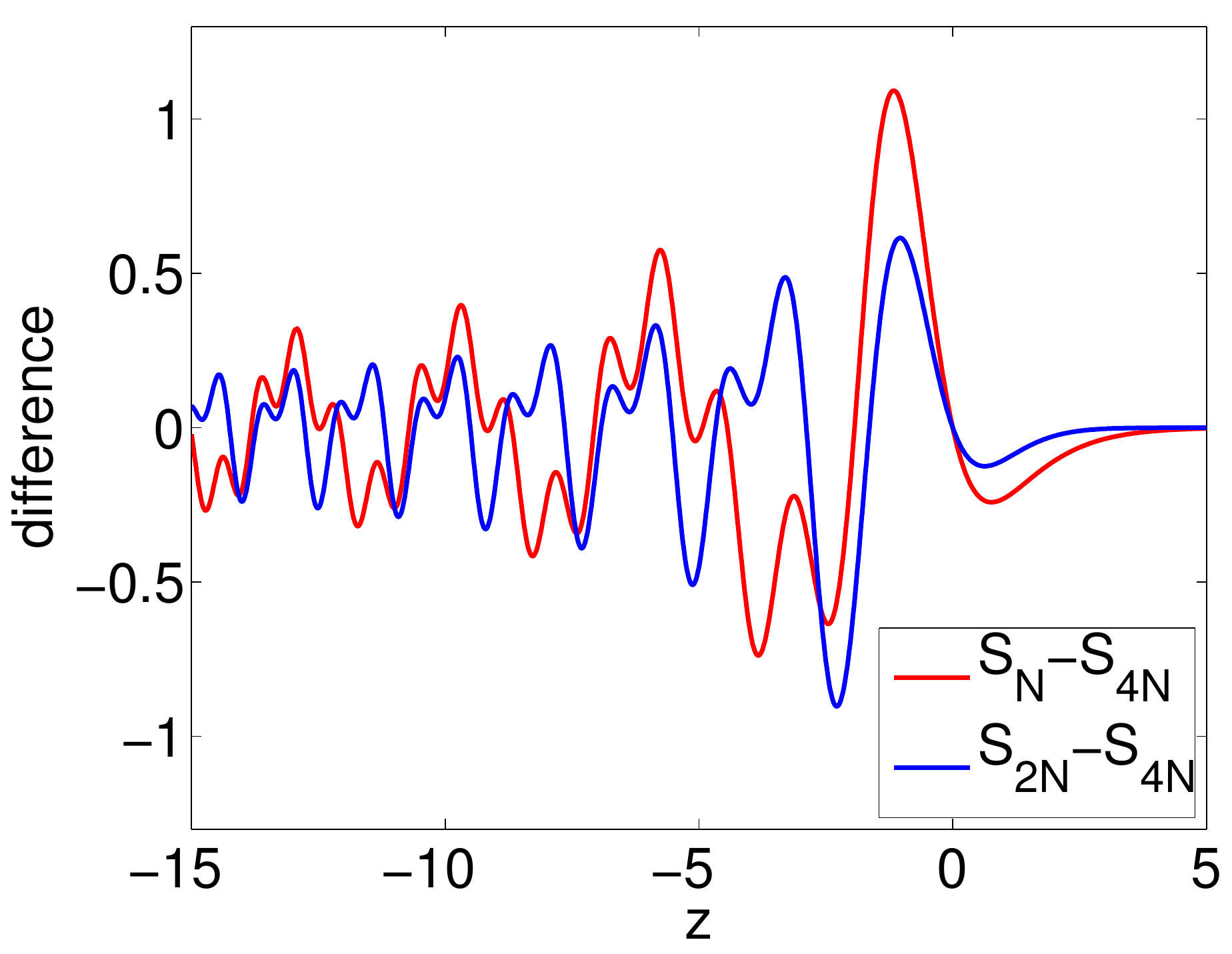}
  \includegraphics[width=2.0in]{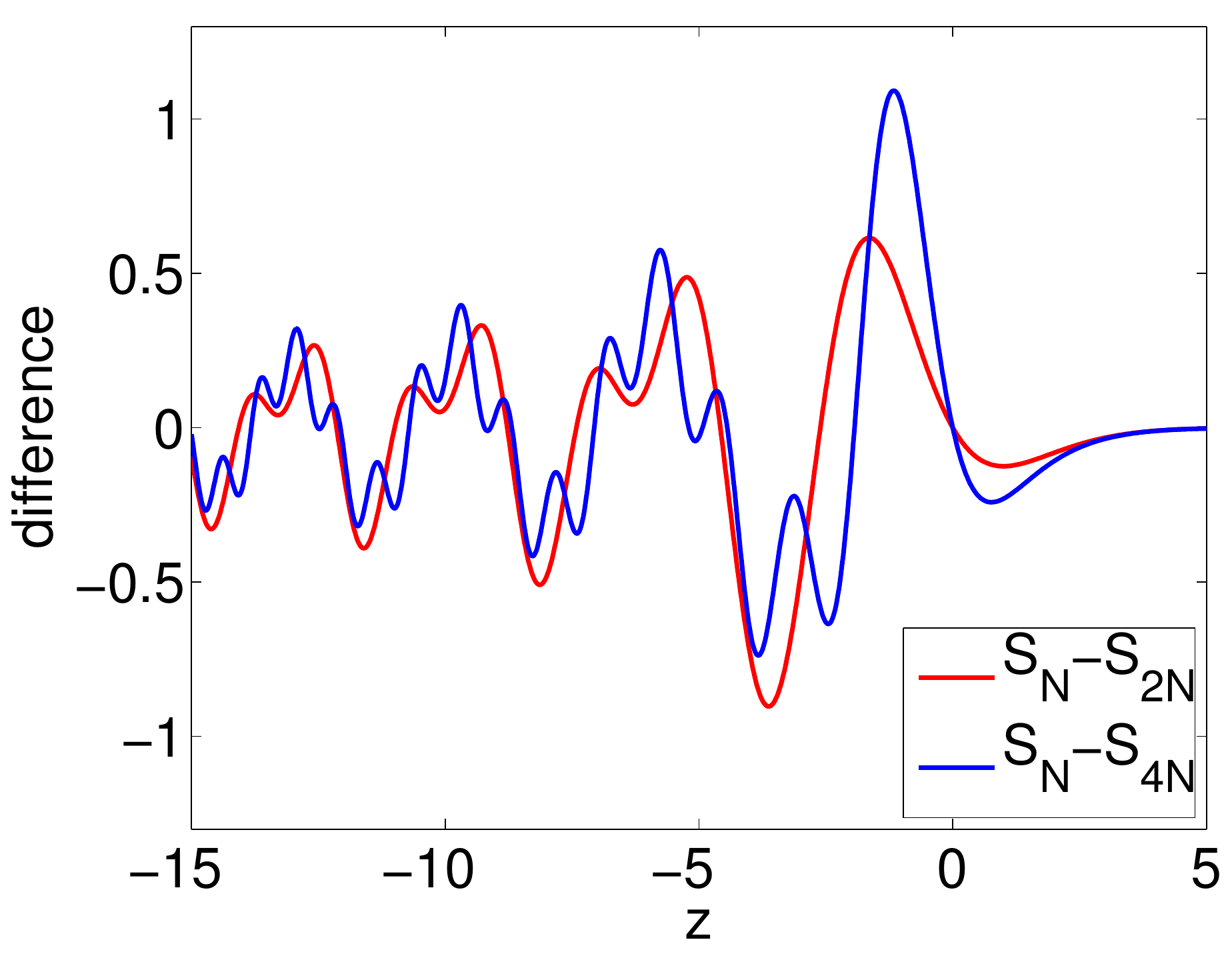}
  \caption{Differences of similarity jumps in the frame of reference moving with the discontinuity.}
  \label{fig:souDiffs}
\end{center}
\end{figure}

\begin{figure}
\begin{center}
  \includegraphics[width=2.0in]{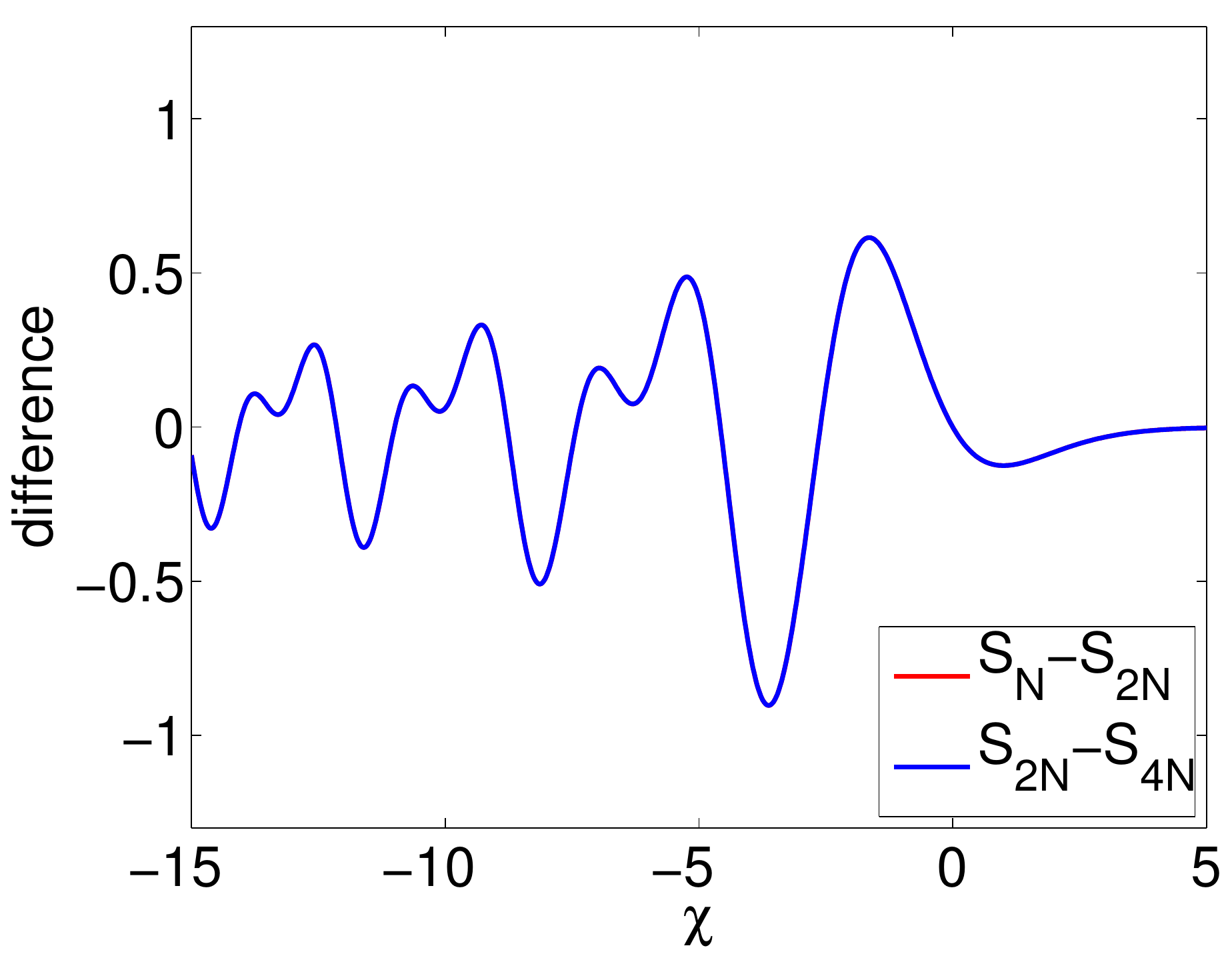}
  \includegraphics[width=2.0in]{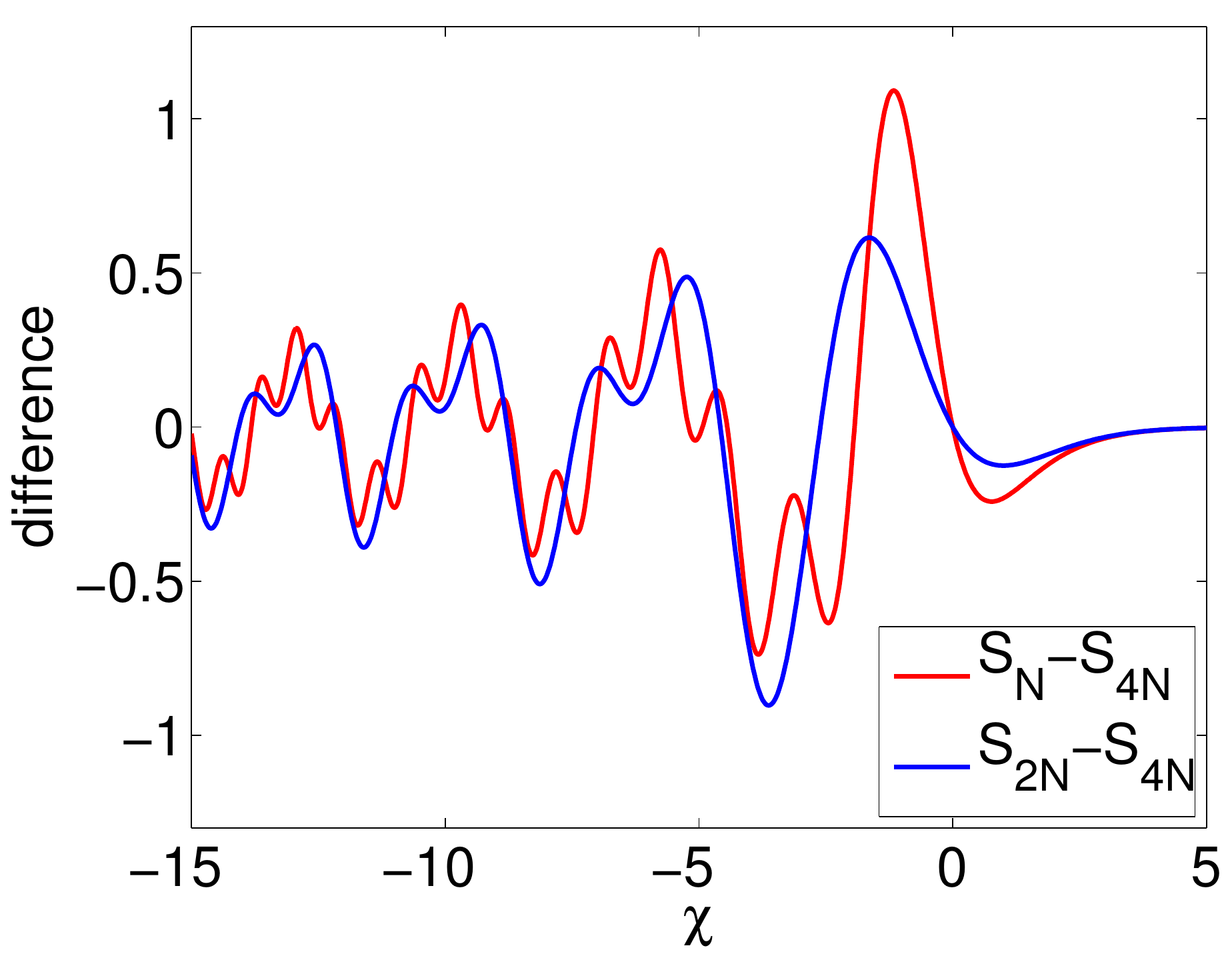}
  \includegraphics[width=2.0in]{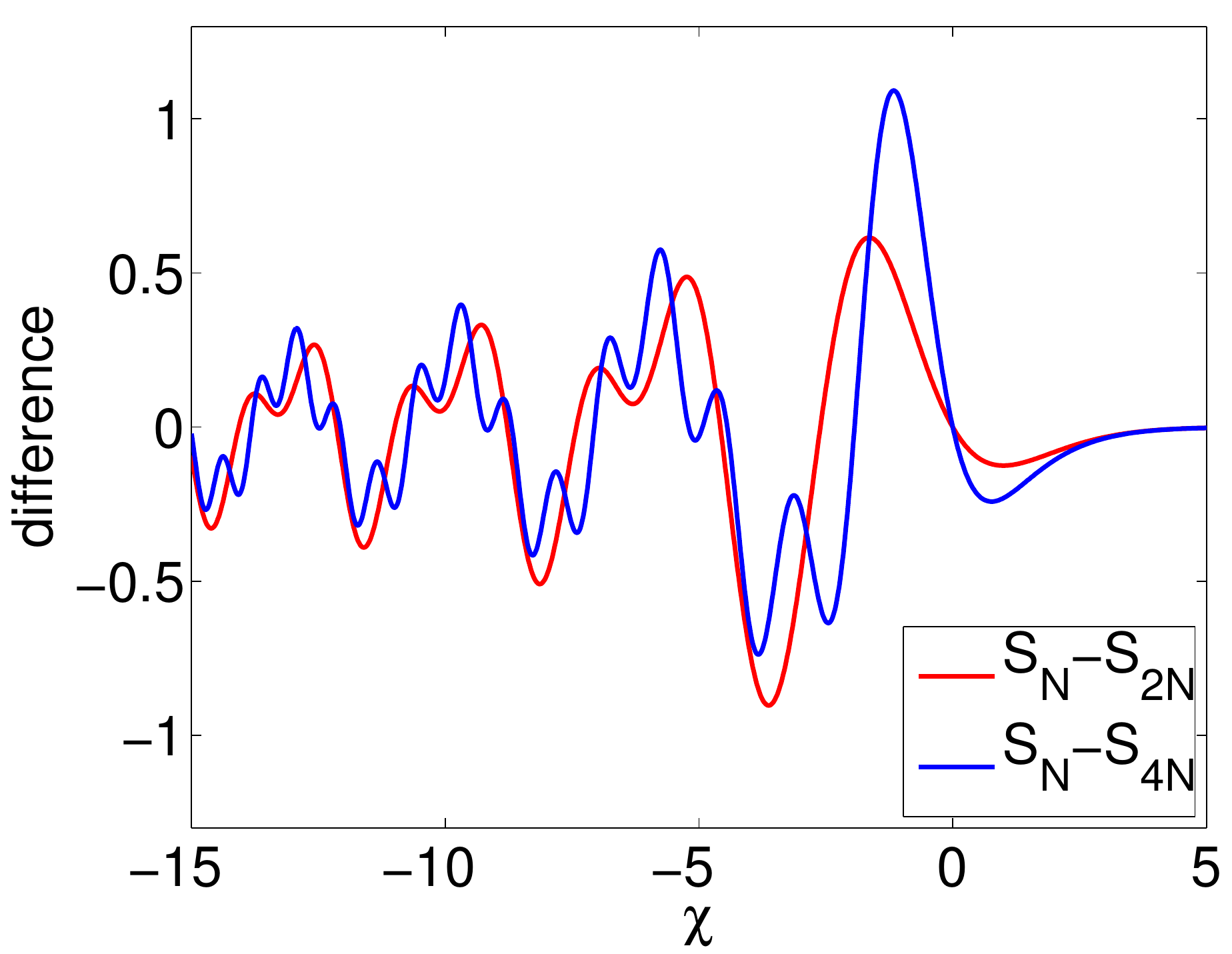}
  \caption{Differences of similarity jumps in the common scaled reference frame $\chi$ (see Equation (\ref{eq:commonChange}) for details). Note that $S_N - S_{2N}$ is nearly identical to $S_{2N} - S_{4N}$ in the upper left plot and is therefore not visible.}
  \label{fig:souDiffs_scaled}
\end{center}
\end{figure}

In order to more clearly explain what is going on, we present additional details for this case. As discussed above in Section~\ref{sec:HO}, the crux of the matter centers around the similarity solution $S$. For this second-order scheme, the solution can be found as
\begin{equation}
  u_h(\xi_h) = \frac{1}{3}-\frac{\xi_h\left(
    \xi_h\sqrt{3}\left(\Gamma \left(\frac{2}{3}\right)\right)^2\, {_1F_2}\left(\frac{2}{3};\frac{4}{3},\frac{5}{3};\frac{\xi_h^3}{27}\right) - 4\pi\, {_1F_2}\left(\frac{1}{3};\frac{2}{3},\frac{4}{3};\frac{\xi_h^3}{27}\right)
    \right)}{6\Gamma \left(\frac{2}{3}\right)\pi}
  \label{eq:HOsoln}
\end{equation}
where $\Gamma$ is the Euler Gamma function, and $_1F_2$ is a generalized hypergeometric function. Figure~\ref{fig:souJump} shows similarity solutions in the reference frame moving with the discontinuity at three resolutions. The grid spacing is essentially a parameter, and so we have chosen a normalization $h_1=1$. Solutions with two grid doublings are also shown. 
Following the analysis in Section~\ref{sec:HO}, differences of the three solutions will be taken. Figure~\ref{fig:souDiffs} shows the three sets of differences which are produced for the three variants of the Richardson extrapolation error estimate. All three plots show the very complex character of the function whose absolute integral is taken. 
The key observation of this paper is presented in Figure~\ref{fig:souDiffs_scaled} where the spatial variable is scaled to the common reference variable $\chi$, as suggested in Equation (\ref{eq:commonChange}). For the case of uniform refinement when the differences are made as suggested, the integrals in the numerator and denominator of (\ref{eq:primaryRatio}) are identical, and the estimate follows. In this case the method essentially avoids the need to calculate the actual integral and relies on the fact that the ratio is known a-priori as one for any similarity function.

\junk{\subsection{Analysis for this case}
Here we end up showing that the convergence rate is dependent on the ratios of resolution and can be negative.

The modified equation for this scheme is
\begin{equation}
  \frac{\partial}{\partial t}\bu+a  \frac{\partial}{\partial x}\bu
    +\frac{a h^2}{2}\left(1-\cfl \right)\left(\frac{1}{3}(1+\cfl)-\frac{1}{2}\right)
    \frac{\partial^3}{\partial x^3}\bu
    +\cdots=0.
    \label{eq:mod2ndUpwind}
\end{equation}
Truncating the series gives
\begin{equation}
  \frac{\partial}{\partial t} \bu_N+a  \frac{\partial}{\partial x} \bu_N+\eta_N \frac{\partial^3}{\partial x^3}\bu_N=0,
    \label{eq:modtrunc2ndUpwind}
\end{equation}
where $\eta_N$ is
\begin{equation}
  \eta_N=\frac{a h_N^2}{2}\left(1-\cfl \right)\left(\frac{1}{3}(1+\cfl)-\frac{1}{2}\right).
\end{equation}
Note that $\eta$ changes sign at $\lambda=\frac{1}{2}$ and scheme becomes third order accurate at $\lambda=\frac{1}{2}$. The case $\lambda=\frac{1}{2}$ is outside our analysis.  For discontinuous initial data $u(x,0)=u_L$ for $x<x_0$ and $u(x,0)=u_R$ for $x  \ge x_0$, one can seek a similarity solution. This approach is detailed in~\cite{banks08a} and yields
\begin{equation}
\begin{array}{l}
  \displaystyle{
    u_N(\xi)=\frac{1}{3}(u_L+2 u_R)} \medskip\\
  \hspace{0.25in} \displaystyle{ -\frac{(u_R-u_L)\xi  \left(\sqrt[3]{\frac{1}{\eta }} \xi  \Gamma \left(\frac{2}{3}\right)^2 \,
   _1F_2\left(\frac{2}{3};\frac{4}{3},\frac{5}{3};\frac{\xi ^3}{27 \eta }\right)-3 \Gamma \left(\frac{1}{3}\right) \Gamma \left(\frac{5}{3}\right) \,
   _1F_2\left(\frac{1}{3};\frac{2}{3},\frac{4}{3};\frac{\xi ^3}{27 \eta }\right)\right)}{27 \sqrt[3]{\eta } \Gamma \left(\frac{2}{3}\right) \Gamma
   \left(\frac{4}{3}\right) \Gamma \left(\frac{5}{3}\right)}}
   \end{array}
   \label{eq:HOsoln}
\end{equation}
where $\Gamma$ is the Euler Gamma function, and $_1F_2$ is a generalized hypergeometric function.  The analytic solution to (\ref{eq:modtrunc2ndUpwind}) can be written by simply substituting $\xi=\frac{x-x_0-a t}{t^{1/3}}$ into (\ref{eq:HOsoln}). In order to simplify the discussion as much as possible, we rewrite Equation~\eqref{eq:HOsoln} in the generic form 
\[
  u_N(\xi) = C_1+C_2S(\xi;\eta_N)
\]
where $C_1$ and $C_2$ depend only on the left and right states, while $S(\xi;\eta_N)$ is the similarity function
\[
  S(\xi;\eta_N) = -\xi \sqrt[3]{\frac{1}{\eta}}  \left(\sqrt[3]{\frac{1}{\eta }} \xi  \Gamma \left(\frac{2}{3}\right)^2 \,
   _1F_2\left(\frac{2}{3};\frac{4}{3},\frac{5}{3};\frac{\xi ^3}{27 \eta }\right)-3 \Gamma \left(\frac{1}{3}\right) \Gamma \left(\frac{5}{3}\right) \,
   _1F_2\left(\frac{1}{3};\frac{2}{3},\frac{4}{3};\frac{\xi ^3}{27 \eta }\right)\right).
\]
Now,
\begin{align*}
  \norm{u_M(x)-u_N(x)}_1 & = \int_{-\infty}^{\infty} \left|C_2S(\xi;\eta_M)-C_2S(\xi,\eta_N)\right| \, d\xi
\end{align*}
from whence
\begin{align*}
  \frac{\norm{u_1(x)-u_3(x)}_1}{\norm{u_2(x)-u_3(x)}_1} & = \frac{\int_{-\infty}^{\infty} \left|S(\xi;\eta_1)-S(\xi,\eta_3)\right| \, d\xi}{\int_{-\infty}^{\infty} \left|S(\xi;\eta_2)-S(\xi,\eta_3)\right| \, d\xi}.
\end{align*}
The difficulty here is that unlike the monotone case, the two functions inside the absolute values oscillate around each other and so no simple choice can be made which removes the absolute value bars. Furthermore, the integrals can be made to converge to many (possibly any) finite number by picking ratios of grid spacing. More on this to come ...}

%-------------- additional 1D examples
\subsection{Additional demonstration of the theory in 1D}
\label{sec:examples}

In order to further demonstrate the validity of the theory we have just described, we perform a series of tests for linear schemes of increasing order, as well as a high-resolution nonlinear TVD scheme. The linear schemes we investigate here are upwind biased single step schemes (i.e. in advancing from $t^n$ to $t^{n+1}$ they use data from $t^n$ only) that are high-order accurate in space and time. Derivation of these schemes follows standard procedures using a Cauchy-Kowalewski process, sometimes called the Lax-Wendroff procedure~\cite{lax60}. Further details on these derivations can be found in~\cite{despres08,despres09,banks12c}. The nonlinear TVD discretization is of the high-resolution Godunov type described in~\cite{vanleer79,sweby84}.

%-------------- fourth-order
\subsubsection{Fourth-Order Linear Scheme}
\label{subsec:fourth}
Consider the linear fourth-order upwind method
\begin{equation}
  u_i^{n+1}=u_i^n+\sum_{s=-3}^2C^{(4)}_{4+s}u_{i+s}^n
  \label{eq:4thUpwind}
\end{equation}
where $C^{(4)}$ is a vector of stencil coefficients given by
\[
C^{(4)}=\frac{\cfl}{144}
\left[
  \begin{array}{c}
    5-8\cfl^2+3\cfl^3 \\
    -37-6\cfl+52\cfl^2-9\cfl^3 \\
    146+96\cfl-104\cfl^2+6\cfl^3 \\
    -50-180\cfl+80\cfl^2+6\cfl^3 \\
    -71+96\cfl-16\cfl^2-9\cfl^3 \\
    7-6\cfl-4\cfl^2+3\cfl^3
  \end{array}
\right].
\]
Table~\ref{table:4thOrderEstimates} shows results using the three basic variants of Richardson extrapolation, and various choices of the uniform refinement ratio $r$ using (\ref{eq:4thUpwind}). The test problem is again defined $a=1$, $[x_L,x_R]=[-\pi,\pi]$, $t_f=2$, and $\cfl=0.6$. The first column shows very good agreement between the estimated convergence rate and the expected rate of $\frac{p}{p+1}=0.8$. The other two columns indicate that results which are difficult to interpret can be obtained for other procedures. Note that the estimated rates are slightly higher than the theoretical prediction. This behavior has been observed in other studies, and is the result of very slow convergence of high-order numerical approximations to the limiting similarity solutions. Despite the large number of grid points used in this study, we attribute the slight overestimation of the convergence rates in the first column to a lack of sufficient grid resolution.
\begin{table}[hbt]
\begin{center}
\begin{tabular}{|c||c|c|c|} \hline\hline
  \mstrut$r$ & \footnotesize{$\RE{u_{h_1}}{u_{rh_1}}{u_{r^2h_1}}$} &\footnotesize{$\RE{u_{h_1}}{u_{r^2h_1}}{u_{rh_1}}$} & \footnotesize{$\RE{u_{rh_1}}{u_{h_1}}{u_{r^2h_1}}$} \\ \hline\hline
  \mstrut$\frac{1}{2}$   & $0.86$ & $0.23$ & $2.32$ \\ \hline
  \mstrut$\frac{2}{5}$   & $0.83$ & $0.41$ & $2.10$ \\ \hline
  \mstrut$\frac{1}{3}$   & $0.83$ & $0.53$ & $1.91$ \\ \hline
  \mstrut$\frac{2}{7}$   & $0.84$ & $0.73$ & $1.60$ \\ \hline
  \mstrut$\frac{1}{4}$   & $0.85$ & $0.65$ & $1.34$ \\ \hline\hline
\end{tabular}
\caption{Estimated convergence rates for the fourth-order upwind scheme for a solution with a discontinuity. The base resolution uses $51201$ points, and uniform refinement is carried out using a ratio of $r$. Results using the three independent variants of Richardson extrapolation are presented in the various columns.}
\label{table:4thOrderEstimates}
\end{center}
\end{table}

%-------------- sixth-order
\subsubsection{Sixth-Order Linear Scheme}
\label{subsec:sixth}
Consider the linear sixth-order upwind method
\begin{equation}
  u_i^{n+1}=u_i^n+\sum_{s=-4}^3C^{(6)}_{5+s}u_{i+s}^n
  \label{eq:6thUpwind}
\end{equation}
where $C^{(6)}$ is a vector of stencil coefficients given by
\[
C^{(6)}=\frac{\cfl}{4320}
\left[
  \begin{array}{c}
    -31+43\cfl^2-15\cfl^4+3\cfl^5 \\
    289+24\cfl-391\cfl^2-30\cfl^3+123\cfl^4-15\cfl^5 \\
    -1299-324\cfl+1623\cfl^2+360\cfl^3-387\cfl^4+27\cfl^5\\
    4325+3240\cfl-2675\cfl^2-1170\cfl^3+615\cfl^4-15\cfl^5 \\
    -1085-5880\cfl+1505\cfl^2+1680\cfl^3-525\cfl^4-15\cfl^5 \\
    -2589+3240\cfl+267\cfl^2-1170\cfl^3+225\cfl^4+27\cfl^5 \\
    431-324\cfl-419\cfl^2+360\cfl^3-33\cfl^4-15\cfl^5 \\
    -41+24\cfl+47\cfl^2-30\cfl^3-3\cfl^4+3\cfl^5
  \end{array}
\right].
\]
Table~\ref{table:6thOrderEstimates} shows results using the three basic variants of Richardson extrapolation, and various choices of the uniform refinement ratio $r$ using (\ref{eq:6thUpwind}). The test problem is again defined using $a=1$, $[x_L,x_R]=[-\pi,\pi]$, $t_f=2$, and $\cfl=0.6$. The first column again shows remarkable agreement between the estimated convergence rate and the expected rate of $\frac{p}{p+1}=0.86$. The other two columns indicate that results which are difficult to interpret can be obtained for other procedures. As in Section~\ref{subsec:fourth}, we attribute the slight overestimation of the convergence rates in the first column to a lack of sufficient grid resolution.
\begin{table}[hbt]
\begin{center}
\begin{tabular}{|c||c|c|c|} \hline\hline
  \mstrut$r$ & \footnotesize{$\RE{u_{h_1}}{u_{rh_1}}{u_{r^2h_1}}$} &\footnotesize{$\RE{u_{h_1}}{u_{r^2h_1}}{u_{rh_1}}$} & \footnotesize{$\RE{u_{rh_1}}{u_{h_1}}{u_{r^2h_1}}$} \\ \hline\hline
  \mstrut$\frac{1}{2}$   & $0.90$ & $0.16$ & $2.95$ \\ \hline
  \mstrut$\frac{2}{5}$   & $0.90$ & $0.47$ & $2.35$ \\ \hline
  \mstrut$\frac{1}{3}$   & $0.88$ & $0.60$ & $1.93$ \\ \hline
  \mstrut$\frac{2}{7}$   & $0.88$ & $0.78$ & $1.20$ \\ \hline
  \mstrut$\frac{1}{4}$   & $0.90$ & $0.83$ & $1.24$ \\ \hline\hline
\end{tabular}
\caption{Estimated convergence rates for the sixth-order upwind scheme for a solution with a discontinuity. The base resolution uses $51201$ points, and uniform refinement is carried out using a ratio of $r$. Results using the three independent variants of Richardson extrapolation are presented in the various columns.}
\label{table:6thOrderEstimates}
\end{center}
\end{table}
%

%-------------- TVD
\subsubsection{Second Order Nonlinear Scheme}
\label{subsec:minmod}

Finally, we consider a high-resolution TVD limited scheme. The scheme is a second-order MUSCL type scheme using a MinMod limiter applied to the slopes.  The scheme can be written
\begin{equation}
  u_i^{n+1}=u_i^n-\cfl \left[
    \left(u_i^n+\frac{1}{2}(1-\cfl)\alpha\right)-
    \left(u_{i-1}^n+\frac{1}{2}(1-\cfl)\beta\right)
  \right]
  \label{eq:2ndUpwindLimited}
\end{equation}
where
\[
  \alpha=\mbox{MinMod}(u_{i+1}^n-u_i^n,u_i^n-u_{i-1}^n),
\]
\[
  \beta=\mbox{MinMod}(u_{i}^n-u_{i-1}^n,u_{i-1}^n-u_{i-2}^n),
\]
and
\[
  \mbox{MinMod}(b,c)=\left\{
  \begin{array}{lcl}
    b \hspace{0.1in} & \hbox{if} & \hspace{0.1in} |b|<|c| \mbox{ and } bc>0 \\
    c \hspace{0.1in} & \hbox{if} & \hspace{0.1in} |b| \ge |c| \mbox{ and } bc>0 \\
    0 \hspace{0.1in} & \hbox{if} & \hspace{0.1in} bc\le0.
  \end{array}
  \right.
\]
Note that this is nothing more than a high-resolution Godunov method (see \cite{vanleer79,sweby84} for details). Table~\ref{table:TVDOrderEstimates} shows the results for the three variants of Richardson extrapolation convergence estimation. The results in~\cite{banks08a} established that although this scheme is nonlinear, the modified equation has solutions of the form (\ref{eq:genSimilarityForm}). Therefore, our analysis still applies to this case, and the method is expected to yield accurate estimates of the convergence rates for uniform refinement case $\RE{u_{h_1}}{u_{rh_1}}{u_{r^2h_1}}$. 
\begin{table}[hbt]
\begin{center}
\begin{tabular}{|c||c|c|c|} \hline\hline
  \mstrut$r$ & \footnotesize{$\RE{u_{h_1}}{u_{rh_1}}{u_{r^2h_1}}$} &\footnotesize{$\RE{u_{h_1}}{u_{r^2h_1}}{u_{rh_1}}$} & \footnotesize{$\RE{u_{rh_1}}{u_{h_1}}{u_{r^2h_1}}$} \\ \hline\hline
  \mstrut$\frac{1}{2}$   & $0.48$ & $0.48$ & $0.48$ \\ \hline
  \mstrut$\frac{2}{5}$   & $0.55$ & $0.55$ & $0.56$ \\ \hline
  \mstrut$\frac{1}{3}$   & $0.57$ & $0.57$ & $0.57$ \\ \hline
  \mstrut$\frac{2}{7}$   & $0.57$ & $0.57$ & $0.57$ \\ \hline
  \mstrut$\frac{1}{4}$   & $0.60$ & $0.59$ & $0.60$ \\ \hline\hline
\end{tabular}
\caption{Estimated convergence rates for the high-resolution TVD scheme for a solution with a discontinuity. The base resolution uses $51201$ points, and uniform refinement is carried out using a ratio of $r$. Results using the three independent variants of Richardson extrapolation are presented in the various columns.}
\label{table:TVDOrderEstimates}
\end{center}
\end{table}
The results in Table~\ref{table:TVDOrderEstimates} show this to be largely true with notable qualifications. Notice that the estimates become more accurate with decreasing $r$. As $r$ becomes smaller, the finest resolution becomes relatively more accurate, and one might expect an extrapolation estimate to yield results which are more similar to those found when comparing computed results to the exact solution. Indeed, as seen from the prior results in Tables~\ref{table:2ndOrderEstimates}, \ref{table:4thOrderEstimates} and~\ref{table:6thOrderEstimates}, even poorly constructed extrapolation techniques tend to yield somewhat more accurate results as $r$ decreases. In addition, there is some significance to the fact that the numerical values in table~\ref{table:TVDOrderEstimates} are somewhat less accurate than the results from the second-order linear scheme in table~\ref{table:2ndOrderEstimates}. The root cause of this observation is not entirely understood, but our conjecture is that the discontinuities in the higher derivatives of the similarity solution for the MinMod scheme~\cite{banks08a} result in slow convergence. Finally, we note that all the estimates in table~\ref{table:TVDOrderEstimates}, including those given by $\RE{u_{h_1}}{u_{r^2h_1}}{u_{rh_1}}$ and $\RE{u_{rh_1}}{u_{h_1}}{u_{r^2h_1}}$, appear equally accurate. This is not in contradiction with the theory which says that the estimate $\RE{u_{h_1}}{u_{rh_1}}{u_{r^2h_1}}$ will be accurate but not that the others will be inaccurate. In fact this is a similar result to those results for the first-order scheme in Table~\ref{table:1stOrderEstimates}. As was the case for the first-order discretization, this fortuitous behavior can be traced to the monotonicity of the approximations.

%-------------- return to sosup topHat
\section{Revisiting the top-hat problem}\label{sec:tophat}
Having demonstrated the behavior of Richardson extrapolation error estimation for linear jumps in 1D, we now return to the introductory example of Section~\ref{sec:sosup}. Because the second-order wave equation can be written as a system of first-order hyperbolic equations, the analysis presented in the previous sections is expected to have applicability here also. We therefore conduct a similar set of studies to those presented in Section~\ref{sec:model}. We perform extrapolation estimation to a series of computations that use a uniform refinement process and present results for the three variants of \Rich while varying the refinement ratio $r$. Tables~\ref{table:TH1st} through~\ref{table:TH4th} present these results.
\begin{table}[hbt]
\begin{center}
\begin{tabular}{|c||c|c|c|} \hline\hline
  \mstrut$r$ & \footnotesize{$\RE{u_{h_1}}{u_{rh_1}}{u_{r^2h_1}}$} &\footnotesize{$\RE{u_{h_1}}{u_{r^2h_1}}{u_{rh_1}}$} & \footnotesize{$\RE{u_{rh_1}}{u_{h_1}}{u_{r^2h_1}}$} \\ \hline\hline
  \mstrut$\frac{1}{2}$   & $0.44$ & $0.40$ & $0.50$ \\ \hline
  \mstrut$\frac{2}{5}$   & $0.47$ & $0.43$ & $0.53$ \\ \hline
  \mstrut$\frac{1}{3}$   & $0.47$ & $0.44$ & $0.53$ \\ \hline
  \mstrut$\frac{2}{7}$   & $0.47$ & $0.44$ & $0.54$ \\ \hline
  \mstrut$\frac{1}{4}$   & $0.47$ & $0.44$ & $0.53$ \\ \hline\hline
\end{tabular}
\caption{Estimated convergence rates for the first-order upwind scheme and the top-hat problem. The base resolution has 400 grid points per dimension and uniform refinement with ratio $r$ is carried out. Results using the three independent variants of Richardson extrapolation are presented in the various columns.}
\label{table:TH1st}
\end{center}
\end{table}
\begin{table}[hbt]
\begin{center}
\begin{tabular}{|c||c|c|c|} \hline\hline
  \mstrut$r$ & \footnotesize{$\RE{u_{h_1}}{u_{rh_1}}{u_{r^2h_1}}$} &\footnotesize{$\RE{u_{h_1}}{u_{r^2h_1}}{u_{rh_1}}$} & \footnotesize{$\RE{u_{rh_1}}{u_{h_1}}{u_{r^2h_1}}$} \\ \hline\hline
  \mstrut$\frac{1}{2}$   & $0.53$ & $-0.27$ & $1.92$ \\ \hline
  \mstrut$\frac{2}{5}$   & $0.56$ & $-0.08$ & $2.02$ \\ \hline
  \mstrut$\frac{1}{3}$   & $0.55$ & $0.06$ & $1.90$ \\ \hline
  \mstrut$\frac{2}{7}$   & $0.57$ & $0.17$ & $1.87$ \\ \hline
  \mstrut$\frac{1}{4}$   & $0.57$ & $0.25$ & $1.73$ \\ \hline\hline
\end{tabular}
\caption{Estimated convergence rates for the second-order upwind scheme and the top-hat problem. The base resolution has 400 grid points per dimension and uniform refinement with ratio $r$ is carried out. Results using the three independent variants of Richardson extrapolation are presented in the various columns.}
\label{table:TH2nd}
\end{center}
\end{table}
\begin{table}[hbt]
\begin{center}
\begin{tabular}{|c||c|c|c|} \hline\hline
  \mstrut$r$ & \footnotesize{$\RE{u_{h_1}}{u_{rh_1}}{u_{r^2h_1}}$} &\footnotesize{$\RE{u_{h_1}}{u_{r^2h_1}}{u_{rh_1}}$} & \footnotesize{$\RE{u_{rh_1}}{u_{h_1}}{u_{r^2h_1}}$} \\ \hline\hline
  \mstrut$\frac{1}{2}$   & $0.62$ & $0.52$ & $0.79$ \\ \hline
  \mstrut$\frac{2}{5}$   & $0.63$ & $0.55$ & $0.78$ \\ \hline
  \mstrut$\frac{1}{3}$   & $0.63$ & $0.57$ & $0.75$ \\ \hline
  \mstrut$\frac{2}{7}$   & $0.63$ & $0.58$ & $0.76$ \\ \hline
  \mstrut$\frac{1}{4}$   & $0.63$ & $0.58$ & $0.74$ \\ \hline\hline
\end{tabular}
\caption{Estimated convergence rates for the nonlinear high-resolution upwind scheme and the top-hat problem. The base resolution has 400 grid points per dimension and uniform refinement with ratio $r$ is carried out. Results using the three independent variants of Richardson extrapolation are presented in the various columns.}
\label{table:THHR}
\end{center}
\end{table}
\begin{table}[hbt]
\begin{center}
\begin{tabular}{|c||c|c|c|} \hline\hline
  \mstrut$r$ & \footnotesize{$\RE{u_{h_1}}{u_{rh_1}}{u_{r^2h_1}}$} &\footnotesize{$\RE{u_{h_1}}{u_{r^2h_1}}{u_{rh_1}}$} & \footnotesize{$\RE{u_{rh_1}}{u_{h_1}}{u_{r^2h_1}}$} \\ \hline\hline
  \mstrut$\frac{1}{2}$   & $0.64$ & $-0.50$ & $3.41$ \\ \hline
  \mstrut$\frac{2}{5}$   & $0.65$ & $-0.11$ & $3.36$ \\ \hline
  \mstrut$\frac{1}{3}$   & $0.65$ & $0.11$ & $2.85$ \\ \hline
  \mstrut$\frac{2}{7}$   & $0.66$ & $0.25$ & $2.63$ \\ \hline
  \mstrut$\frac{1}{4}$   & $0.66$ & $0.34$ & $2.33$ \\ \hline\hline
\end{tabular}
\caption{Estimated convergence rates for the fourth-order upwind scheme and the top-hat problem. The base resolution has 400 grid points per dimension and uniform refinement with ratio $r$ is carried out. Results using the three independent variants of Richardson extrapolation are presented in the various columns.}
\label{table:TH4th}
\end{center}
\end{table}

The results show that the conclusions derived from the simple 1D analysis hold even for this more complex situation. In particular, the estimated rate of convergence found when using uniform refinement and using $\RE{u_{h_1}}{u_{rh_1}}{u_{r^2h_1}}$ yields results that are in reasonable agreement with the a-priori expected rate $p/(p+1)$. However, the two other variants are found to yield wildly varying results for both the second-order and fourth-order cases. As for the 1D case, the first-order and high-resolution methods yield monotone (or nearly monotone) approximations, which leads to fortuitous cancelation. As a result the estimated convergence rates are nearly correct for all three methods. This behavior is not a general result, and should not be relied upon in practice. One additional interesting point is that as the separation between the grid spacing of the three numerical solutions increases, the estimated convergence rate seems to converge toward the a-priori expectation. For example, the estimates in the second and third columns in Table~\ref{table:TH4th} are woefully inaccurate, but in both cases the accuracy is improving as the relative separation increases. This may be one reason why practical studies, such as those in~\cite{appelo12}, have used \Rich with extremely refined final simulations.

%-------------- Conclusions
\section{Conclusions}
\label{sec:conclusions}
We have provided an in depth investigation of Richardson extrapolation error estimation for linear hyperbolic wave propagation in the presence of discontinuous solutions. The discussion is motivated by an example of 2D scalar wave propagation that has been taken from the literature. Numerical methods of orders 1,2,4, and a high-resolution (nominally $2^{\hbox{nd}}$ order) variant are used to produce approximations at three resolutions. Richardson extrapolation error estimates found using these approximations are found to produce results that can vary widely from the expected convergence rate of $p/(p+1)$ where $p$ is the order of the method. One exception to this wide variation is observed if the estimate is determined in a specific manner.

In order to understand this behavior, a simpler model problem of 1D advection is posed as an appropriate surrogate, and a detailed analysis is presented. The analysis uses the solution to the modified equation to elucidate the difficulty found in practice. In addition, the analysis reveals that one particular realization of the technique reproduces the a-priori convergence rates even in the presence of a discontinuity or other similarity type behavior. The key elements are shown to be the use of uniform refinement, and that the comparisons inherent in the Richardson estimator are performed in one specific manner as $\RE{u_{h_1}}{u_{rh_1}}{u_{r^2h_1}}$ where $r$ is the uniform refinement ratio. This result was then demonstrated in 1D for a number of discretizations ranging in order from first-order to sixth-order, and for the motivating 2D discretizations ranging from first-order to fourth-order. In addition, results were presented for nonlinear high-resolution schemes. The results were found to be in good agreement with the theory. 

The analysis presented in this work is an interesting result that shows why estimated convergence rates from Richardson extrapolation for problems with discontinuities have been found to be difficult to interpret. In addition, it shows that Richardson extrapolation can be used to obtain predictable a-priori convergence rates even for simulations of systems with linear discontinuities provided two key elements are satisfied. First, uniform refinement must be used. Second, the estimate must be obtained as $\RE{u_{h_1}}{u_{rh_1}}{u_{r^2h_1}}$ in the notation of this paper. The result is in fact more general than discontinuities and includes other self-similar features such as corners. Future work will include investigating the techniques discussed here for nonlinear equations such as the Euler equations, as well parabolic problems or dispersive wave propagation problems.

\bibliographystyle{spphys}
\bibliography{journal-ISI,jwb}

\end{document}